\documentclass[10pt]{amsart}
\calclayout
\usepackage[T1]{fontenc}
\usepackage{lmodern}
\usepackage{amsmath,amssymb,mathtools}
\usepackage{microtype}
\usepackage{tikz}
\usepackage[hidelinks]{hyperref}
\numberwithin{equation}{section}
\newtheorem{theorem}{Theorem}[section]
\newtheorem{lemma}[theorem]{Lemma}
\newtheorem{proposition}[theorem]{Proposition}
\newtheorem{corollary}[theorem]{Corollary}

\theoremstyle{definition}

\newtheorem{question}[theorem]{Question}
\theoremstyle{remark}

\newcommand{\N}{\omega}
\newcommand{\Cantor}{2^\omega}
\newcommand{\Strings}{2^{<\omega}}
\newcommand{\restr}{\mathbin{\upharpoonright}}
\newcommand{\supp}{\operatorname{supp}}
\newcommand{\dom}{\operatorname{dom}}
\newcommand{\Prob}{\mathbb{P}}
\DeclareMathOperator{\dimH}{dim}
\title[Effective dimension one and the $K$-degrees]{Effective Hausdorff dimension one is cofinal in the $K$-degrees}
\author{Xuanheng Zhao}
\address{School of Mathematics, Nanjing University,
Nanjing, Jiangsu 210093, People's Republic of China}
\email{xuanheng21@gmail.com}

\subjclass[2020]{Primary 03D32; Secondary 68Q30}
\keywords{Kolmogorov complexity, $K$-degree, effective Hausdorff dimension}
\hypersetup{pdftitle={Effective Hausdorff dimension one is cofinal in the K-degrees}}
\begin{document}
\begin{abstract}
Every $K$-degree has an upper bound of effective Hausdorff
dimension one and Kurtz random
sequences are cofinal in the $K$-degrees.
\end{abstract}
\maketitle
\setcounter{tocdepth}{2}
\begingroup
\small
\endgroup

\section{Introduction}

Algorithmic randomness studies the extent to which an infinite sequence
can resist effective tests and compression. Prefix-free Kolmogorov
complexity gives a way to examine this question at each finite length:
$K(\sigma)$ is the length of a shortest self-delimiting program producing
the binary string $\sigma$, relative to a fixed optimal prefix-free
machine. The function $n\mapsto K(A\restr n)$ records the information
content of the successive initial segments of $A$. Comparing these
functions leads to a finer analysis than merely asking whether a
sequence is random.

For $A,B\in\Cantor$, write
\[
 A\leq_K B\quad\Longleftrightarrow\quad
 (\exists c)(\forall n)\ K(A\restr n)\leq K(B\restr n)+c.
\]
Thus every initial segment of $A$ is no more complex than the initial
segment of $B$ of the same length, apart from one fixed additive
constant. The equivalence relation $A\equiv_K B$ means that both
$A\leq_K B$ and $B\leq_K A$ hold. Its equivalence classes are the
\emph{$K$-degrees}, partially ordered by $\leq_K$. We write $A<_K B$
when $A\leq_K B$ but $B\nleq_K A$. This approach to relative complexity,
with antecedents in Solovay's work, was developed systematically by
Downey, Hirschfeldt and LaForte~\cite{DHL}. Barmpalias~\cite{Barmpalias}
surveys the resulting degree structures and their connections with
recursion theory and randomness.

The least $K$-degree consists of the \emph{$K$-trivial} sequences,
characterized by $K(A\restr n)\leq K(n)+O(1)$. Beyond this least degree,
even the random sequences have a rich order structure. Yu and
Ding~\cite{YD} proved that $1$-random sequences realize
$2^{\aleph_0}$ distinct $K$-degrees. Yu, Ding and
Downey~\cite{YDD} showed that there is no greatest $K$-degree.
Csima and Montalb\'an~\cite{CM} constructed a minimal pair of nonzero
$K$-degrees: their only common lower bound is the least degree.
Miller and Yu~\cite{MY08} established a different obstruction to simple
order structure: if $X\oplus Y$ is $1$-random, then $X$ and
$Y$ have no common $K$-upper bound. Their subsequent work~\cite[Theorem~6.1]{MY11}
showed that the $K$-degrees of $1$-random sequences have
no minimal elements: every random sequence has a random sequence
strictly below it in the $K$-degrees. 
Miller~\cite{Miller09} proved
that the upper cone above any $2$-random sequence is countable.

A central question is whether 
\begin{equation}\label{eq:random-upper-question}
 (\forall A\in\Cantor)(\exists R\in\Cantor)
 \bigl[R\text{ is }1\text{-random and }A\leq_K R\bigr].
\end{equation}
Miller and Nies record it as Question~9.8 in~\cite{MN}.
It is also the $\leq_K$ part of Question~6 in Barmpalias's
survey~\cite{Barmpalias} and appears as Open Question~6.3 in
Hirschfeldt's survey~\cite{Hirschfeldt}.
Downey, Hirschfeldt, Nies and Terwijn~\cite[Section~12.1]{DHNT}
also discuss it.
The present paper gives partial results towards this question.
We do not settle~\eqref{eq:random-upper-question}.

A closely related structural question concerns maximality. A degree
$\mathbf a$ is \emph{maximal} if there is no degree strictly above it.
This differs from being a \emph{greatest} degree, which would have to
bound every degree. Miller and
Nies~\cite[Question~9.7]{MN} asked whether maximal $K$-degrees exist.
Liu~\cite{Liu} recently answered this question negatively:
$(\forall A\in\Cantor)(\exists B\in\Cantor)\ A<_K B$.
However, for a nonrandom source, Liu's theorem does not prescribe a random
representative for the upper degree. The block transformations
in his proof are recursive measure-preserving bijections with recursive
inverses, so they preserve both $1$-randomness and
nonrandomness~\cite[Section~2.2]{Liu}. Consequently, whenever such a
transformation raises the $K$-degree of a nonrandom input, the resulting
upper degree has no random representative, since randomness is invariant
under $\equiv_K$. 
Our results impose two other requirements on upper bounds: effective
Hausdorff dimension one and Kurtz randomness.

Effective Hausdorff dimension measures the asymptotic lower rate of
initial-segment complexity. By the complexity characterization of
Mayordomo~\cite{Mayordomo}, in the theory developed by
Lutz~\cite{Lutz},
$\dimH(X)=\liminf_{n\to\infty}\frac{K(X\restr n)}n$.
Its largest possible value is $1$. Every $1$-random sequence
has dimension $1$, but dimension $1$ allows an unbounded sublinear
complexity deficiency. Our main result is the following.

\begin{theorem}\label{thm:main}
For every $A\in\Cantor$, there is $B\in\Cantor$ such that
$A\leq_K B$, $A\leq_T B$, and $\dimH(B)=1$.
\end{theorem}

Our second result concerns \emph{Kurtz randomness}: a sequence is
Kurtz random if it belongs to no effectively closed null class.
We prove that Kurtz random sequences are cofinal in the $K$-degrees,
and that the comparison constant can be chosen independently of the
source. In fact, one fixed closed class of positive measure consisting
entirely of Kurtz random sequences contains a $K$-upper bound for
every source (Theorem~\ref{kz:cofinal}). We also show that some
$\Delta^0_2$ $K$-degree contains no Kurtz random representative
(Proposition~\ref{kz:no-random-representative}), using Tadaki's exact
complexity theorem~\cite{Tadaki}.

Dimension one has further consequences
for the statistical and computational properties of upper bounds.
Every $K$-degree has an upper bound satisfying the strong law of large
numbers (Corollary~\ref{cor:strong-law}). The sets
of diagonally nonrecursive (DNR) Turing degree are also cofinal in the
$K$-degrees (Corollary~\ref{cor:dnc}). 

We also obtain a uniform-distribution consequence. Every $A$ has a
$K$-upper bound $B$ such that, for the real $x_B\in[0,1]$ with binary
expansion $B$, the sequence $(a_nx_B)_n$ is uniformly distributed
modulo one for every recursive sequence $(a_n)$ of pairwise distinct
integers satisfying $|a_n|=O(n^d)$ for some integer $d\geq1$
(Corollary~\ref{cor:uniform-distribution}), which is a restricted form of Avigad's UD randomness~\cite{Avigad}.

Two further results clarify the original $1$-random upper-bound
question. First, all $K$-degrees have $1$-random upper bounds if and
only if some weakly $2$-generic sequence has one. If this cofinality statement
fails, a counterexample can be chosen to be a low $1$-generic sequence
(Proposition~\ref{prop:generic-reduction}). Second, the upper bounds
cannot in general be required to be weakly Demuth random. We construct
a non-$1$-random sequence $A$ of effective dimension one such that
$A<_K\Omega$ but $A$ has no weakly Demuth random $K$-upper bound
(Proposition~\ref{prop:weak-demuth-obstruction}).

The two main cofinality proofs use a common interval code. For a suitable probability
measure, the code places an output in its support while allowing a
bounded list containing the source prefix to be recovered from the
output prefix of the same length. The Kurtz argument applies this
code to a measure computed from a fixed nonrecursive r.e.\ set
that is low for $K$. Its settling-time function permits a direct
construction of a positive-measure closed class consisting entirely
of Kurtz random sequences. Lowness for $K$ then removes the oracle
from the complexity comparison.

For dimension one, we alternate ordinary coding with repairs. A
repair stores a shortest conditional program for a longer source
prefix and fills the remaining positions with free bits. Even before
the program is complete, its prefixes contribute complexity. On each
finite block, we consider source complexity conditional on the old
output prefix and its complexity. We begin a repair where this
conditional complexity, after subtraction of a suitable linear
function, is maximal. This choice both controls the length of the
repair program and gives a direct lower bound on output complexity
growth. A finite reserve of source information handles the beginning
of each repair, preserving domination at every length. The ordinary
interval code has bounded-list decoding in both directions. Thus,
if repairs stopped, the output would have the same dimension as the
source, contradicting the block growth estimate. Repairs therefore
continue, and their complexity growth rates tend to one. Each repair
also allows an exact source prefix to be recovered from $B$. These
prefixes have unbounded lengths, giving $A\leq_T B$.

Subsection~\ref{sec:prelim} gives the preliminary complexity facts and
the common coding lemma. Section~\ref{sec:kurtz} proves Kurtz
cofinality and gives the example with no Kurtz random representative.
Section~\ref{sec:construction} contains the construction and
verification of Theorem~\ref{thm:main}. Section~\ref{sec:question}
develops the consequences of dimension one and proves the results on
genericity and weak Demuth randomness. It ends with a brief discussion
of analogous questions for plain Kolmogorov complexity.

\subsection{Preliminaries}\label{sec:prelim}

Write $\Strings$ and $\Cantor$ for the finite and infinite binary
sequences, respectively. We write $X_n=X\restr n$ whenever $n\leq |X|$,
and define the subword $X[i:j)$ by $X[i:j)(r)=X(i+r)$ for $0\leq r<j-i$.
Concatenation is denoted by juxtaposition. We write $\sigma\preceq\tau$
when $\sigma$ is a prefix of $\tau$, allowing equality.
For a finite word $\sigma$ and an infinite sequence $X$, put
\[
 0.\sigma=\sum_{i<|\sigma|}\sigma(i)2^{-i-1},\qquad
 x_X=\sum_{i\geq0}X(i)2^{-i-1},\qquad
 I_\sigma=[0.\sigma,0.\sigma+2^{-|\sigma|}].
\]
Thus $x_X\in I_{X_n}$ for every $n$, including either binary expansion
of a dyadic real. For $\sigma\in\Strings$,
$[\sigma]=\{X\in\Cantor:\sigma\prec X\}$, and $\lambda$ denotes
fair-coin measure. Throughout, we write $\log x=\lceil\log_2 x\rceil$ for $x>0$,
where $\log_2$ denotes the usual real-valued binary logarithm. An effectively closed class is a $\Pi^0_1$ subset of
$\Cantor$, equivalently the path set of a recursive binary tree.
A probability measure $\nu$ is \emph{recursive} if the real numbers
$\nu[\sigma]$ are recursive uniformly in $\sigma$. We use the same
definition relative to an oracle $W$.
The \emph{support} of a Borel probability measure $\nu$ on $\Cantor$
is the closed set
\[
 \supp(\nu)=\{X\in\Cantor:(\forall n)\ \nu[X_n]>0\}.
\]
Equivalently, it is the smallest closed set of full $\nu$-measure.

Fix optimal prefix-free machines for unconditional and conditional
complexity. For an oracle $W$, let $K^W$ denote relativized prefix-free
complexity. In particular, $K^W(x)\leq K(x)+O(1)$ by ignoring the oracle.
We reserve $C$ for plain Kolmogorov complexity, with respect to fixed
optimal machines for unconditional and conditional plain complexity.
Finite tuples are identified with their images under a fixed
recursive coding. The constants implicit in $O(1)$ are uniform in the
finite objects involved. When a constant depends on a parameter, this
dependence will be stated.
We recall the facts about prefix-free complexity that will be used below;
see \cite{DH,Nies}. Their relativized forms will also be used.
The chain rule, in conditional form, is
\begin{equation}\label{eq:chain}
 K(x,y\mid z)=K(x\mid z)+K(y\mid z,x,K(x\mid z))+O(1).
\end{equation}
The constant is independent of $x,y,z$, since the coding and decoding
procedures in the usual proof are uniform in the auxiliary input $z$.
Taking $z$ to be the empty string gives the unconditional form.
For each $x$, fix a shortest program $x^*$ by taking the first program
of length $K(x)$ that outputs $x$.
Then $x^*$ and $(x,K(x))$ are recursively interconvertible, so the
condition $(x,K(x))$ in the unconditional form of \eqref{eq:chain}
may be replaced by $x^*$.
We also use
\[
 K(n)\leq 2\log(n+2)+O(1),\qquad
 K(x)\leq |x|+K(|x|)+O(1).
\]

A \emph{discrete semimeasure} is a function
$q:\Strings\to[0,\infty)$ such that $\sum_x q(x)\leq1$.
It is \emph{left-r.e.} if there is a recursive rational-valued function
$(s,x)\mapsto q_s(x)$ such that $q_s(x)\leq q_{s+1}(x)$ and $q(x)=\lim_s q_s(x)$ for all $x\in\Strings$.
The coding theorem~\cite[Theorem~2.2.25, p.~91]{Nies} implies that
every left-r.e.\ discrete semimeasure satisfies
\begin{equation}\label{eq:semimeasure}
 q(x)\leq \kappa_q2^{-K(x)}
\end{equation}
for a constant $\kappa_q$ depending only on $q$.
We will also use Kraft's inequality in the form
$\sum_x2^{-K(x)}\leq1$.

We use the complexity characterization
$\dimH(X)=\liminf_n K(X_n)/n$ of effective Hausdorff dimension
\cite{Lutz,Mayordomo}.
A sequence $X\in\Cantor$ is \emph{$1$-random} if there is a constant
$d\in\N$ such that $K(X_n)\geq n-d$ for every $n \in \N$.
By the Levin--Schnorr theorem, this is equivalent to Martin--L\"of
randomness; see \cite{DH,Nies}.

Finally, suppose that a fixed partial recursive procedure, on input
$y$, outputs a list of at most $d$ objects containing $x$, where $d$ is
fixed. A shortest description of $y$, together with
$\log d$ bits selecting an entry, gives
\begin{equation}\label{eq:boundedlist}
 K(x)\leq K(y)+O(1),\qquad K(x,y)\leq K(y)+O(1).
\end{equation}

The next lemma allows us to encode a sequence in a prescribed closed
set without increasing the length at which its prefixes can be recovered.

\begin{lemma}\label{lem:listcode}
Let $\nu$ be a recursive probability measure on $\Cantor$ such that $\nu[\tau]\leq2\cdot2^{-|\tau|}$ for every $\tau$.
For every $X\in\Cantor$ there is $Y\in\operatorname{supp}(\nu)$
such that a fixed recursive procedure, given $Y_t$, produces a list $\mathcal L(Y_t)$
of at most 4 words containing $X_t$. The procedure is uniform in
an index for $\nu$.
\end{lemma}

The entire statement relativizes to any oracle
$W$.
\begin{proof}
Put $x=x_X$. For each $\tau\in2^t$, define
\[
 \ell_\tau=\sum_{\substack{\rho\in2^t\\\rho<_{\mathrm{lex}}\tau}}\nu[\rho],
 \qquad r_\tau=\ell_\tau+\nu[\tau],\qquad
 J_\tau=[\ell_\tau,r_\tau].
\]
If $\tau_0<_{\mathrm{lex}}\cdots<_{\mathrm{lex}}\tau_{2^t-1}$
enumerates $2^t$, then $\ell_{\tau_0}=0$, $r_{\tau_i}=\ell_{\tau_{i+1}}$ for $i<2^t-1$, and $r_{\tau_{2^t-1}}=1$.
In particular, $\bigcup_{\tau\in2^t}J_\tau=[0,1]$ and
$|J_\tau|=\nu[\tau]$; the interiors of distinct intervals are disjoint.

Starting with $Y_0$ equal to the empty word, choose inductively $Y_{t+1}=Y_tb_t$ by 
\[
 b_t=\min\{b\in\{0,1\}:x\in J_{Y_tb}\text{ and }\nu[Y_tb]>0\}.
\]
The choice is always possible,
and $Y=\bigcup_tY_t\in \supp(\nu)$.

Given $\tau\in2^t$, put $\varepsilon=2^{-t-4}$ and
compute rational $\widehat\ell_\tau,\widehat r_\tau$ such that $|\widehat\ell_\tau-\ell_\tau|\leq\varepsilon$ and $|\widehat r_\tau-r_\tau|\leq\varepsilon$.
Set $\widetilde J_\tau=
 [\widehat\ell_\tau-\varepsilon,\widehat r_\tau+\varepsilon]$.
Then $J_\tau\subseteq\widetilde J_\tau$ and
$|\widetilde J_\tau|\leq\nu[\tau]+4\varepsilon
 \leq\left(2+\frac14\right)2^{-t}$.
Return the recursive list $\mathcal L(\tau)=\{\rho\in2^t:I_\rho\cap\widetilde J_\tau
 \neq\varnothing\}$.
An interval of length at most $(2+1/4)2^{-t}$ meets at most 4
closed intervals of length $2^{-t}$, so $|\mathcal L(\tau)|\leq4$. 
For $\tau=Y_t$, $X_t\in\mathcal L(Y_t)$.
\end{proof}

\section{Kurtz random upper bounds}\label{sec:kurtz}

The interval code from Lemma~\ref{lem:listcode} gives upper bounds
with a weak randomness property. The proof of the following theorem is based on the proof of Brendle et al. \cite[Theorem~5.2]{BBNN}.

\begin{theorem}\label{kz:cofinal}
There are a constant $c$ and a closed class $P$ of measure at least
$2/3$, all of whose members are Kurtz random, such that for every
$Z\in\Cantor$ there is $R\in P$ satisfying $K(Z_n)\leq K(R_n)+c$ for every $n$.
\end{theorem}

\begin{proof}
Fix a nonrecursive r.e.\ set $W$ which is low for $K$: $K(\sigma)\leq K^W(\sigma)+O(1)$ for all $\sigma$.
Such a set exists by the existence of a nonrecursive r.e.\ $K$-trivial
set~\cite{DHNS} and the equivalence of $K$-triviality and lowness for
$K$~\cite{Nies05}.
Let $(W_s)_s$ be a recursive increasing sequence of finite sets with
union $W$, and put
$m(n)=\min\{s:W_s\restr n=W\restr n\}$.
Then $m\leq_T W$, and no total recursive function eventually dominates
$m$.

Let $(Q_e)_e$ enumerate all effectively closed classes, with
$Q_e=\bigcap_s Q_{e,s}$, where the $Q_{e,s}$ are uniformly recursive
clopen sets decreasing in $s$. Set $\varepsilon_n=\frac{2^{-n}}{6(n+1)}$,
\[
 U_n=\bigcup\left\{
 Q_{e,m(n)}:e\leq n,\ 
 \lambda(Q_{e,m(n)})\leq\varepsilon_n
 \right\},
\]
and $U=\bigcup_n U_n$.
The sets $U_n$ are uniformly $W$-recursive clopen sets, $\lambda(U_n)\leq(n+1)\varepsilon_n=\frac{2^{-n}}6$, and $\lambda(U)\leq\frac13$.
If $\lambda(Q_e)=0$, then
\[
 g_e(n)=\min\{s:\lambda(Q_{e,s})\leq\varepsilon_n\}
\]
is a total recursive function. Choose $n\geq e$ with $m(n)>g_e(n)$.
The decreasing approximation gives
$Q_e\subseteq Q_{e,m(n)}\subseteq U_n\subseteq U$.
Thus $U$ contains every effectively closed null class, and every member
of the $\Pi^0_1(W)$ class $P=\Cantor\setminus U$ is Kurtz random.
Moreover, $\lambda(P)\geq2/3$.

For every $\sigma$, put $u_N(\sigma)=\lambda\!\left([\sigma]\cap\bigcup_{n\leq N}U_n\right)$.
These rational numbers are uniformly $W$-recursive, and
$0\leq\lambda([\sigma]\cap U)-u_N(\sigma)
 \leq\sum_{n>N}\frac{2^{-n}}6=\frac{2^{-N}}6$.
Consequently $\lambda(P\cap[\sigma])$ is uniformly $W$-recursive.
Conditioning $\lambda$ on $P$ gives a $W$-recursive probability measure $\nu(E)=\frac{\lambda(E\cap P)}{\lambda(P)}$ with $\supp(\nu)\subseteq P$, and $\nu[\sigma]\leq\frac32\,2^{-|\sigma|}$ for all $\sigma$. 
The relativization of Lemma~\ref{lem:listcode} now gives, for each
$Z\in\Cantor$, an $R\in\supp(\nu)$ such that a fixed $W$-recursive
procedure takes $R_n$ to a bounded list containing $Z_n$. Hence
\[
 K(Z_n)\leq K^W(Z_n)+O(1)
 \leq K^W(R_n)+O(1)
 \leq K(R_n)+O(1)
\]
using lowness for $K$. All auxiliary objects were
fixed before choosing $Z$, so the final constant is independent of
$Z$ and $n$.
\end{proof}

On the other hand, not every $K$-degree contains a Kurtz random
sequence.

\begin{proposition}\label{kz:no-random-representative}
There is a $\Delta^0_2$ sequence $A$ such that no sequence
$K$-equivalent to $A$ is Kurtz random.
\end{proposition}

\begin{proof}
Fix an optimal prefix-free machine $U$, and let $A$ be the binary
expansion of the generalized halting probability
$Z_U(1/2)=\sum_{p\in\dom U}2^{-2|p|}$.
Tadaki's exact complexity theorem \cite[Theorem~18]{Tadaki} gives $K(A_n)=\frac n2+O(1)$.
Since $Z_U(1/2)$ is left-r.e., $A$ is $\Delta^0_2$.
Suppose that $B\equiv_K A$. Then there is a constant $c$ such that $\frac n2-c\leq K(B_n)\leq\frac n2+c$.
There is also a constant $d$, independent of $\sigma$ and $k$, such that $K(\sigma0^k)\leq K(\sigma)+K(k)+d$.
Choose $k\geq1$ so large that $\frac k2>K(k)+2c+d$.
If $B$ contained $k$ consecutive zeros starting at position $n$, then
$B_{n+k}=B_n0^k$, and hence
\[
 \frac{n+k}{2}-c
 \leq K(B_{n+k})
 \leq \frac n2+c+K(k)+d,
\]
contrary to the choice of $k$.
Thus $B$ belongs to the effectively closed null class
\[
 P_k=\{X\in\Cantor:(\forall n)\ X[n:n+k)\neq0^k\}.
\]
Therefore $B$ is not Kurtz random.
\end{proof}

\section{Dimension-one upper bounds}\label{sec:construction}

Fix $A\in\Cantor$. If $\dimH(A)=1$, take $B=A$ and Theorem \ref{thm:main} follows. Henceforth suppose that $d=\dimH(A)<1$, and choose an integer $k_*\geq2$ such that
$1-\frac1{k_*}>d$.

\subsection{Idea}

The construction combines ordinary coding with repairs. Ordinary
coding transfers information from $A$ to $B$ through the interval
code of Lemma~\ref{lem:listcode} for a suitable recursive probability
measure. When the source's conditional complexity grows slowly enough over a
sufficiently long part of a block, we can encode the required source
prefix by a short conditional program. A repair stores this program
and uses the remaining free positions to raise the complexity of $B$.
We arrange these two
parts so that $K(A_n)\leq K(B_n)+O(1)$ at every length and
$\dimH(B)=1$.

At stage $j$, we extend $B_N$ to $B_M$,
where $N=2j^2$ and $M=2(j+1)^2$. We call $[N,M]$ block $j$.
It consists of ordinary coding, possibly followed by one repair on
$[a,M)$ for some $a\in[N,M)$.
Ordinary coding may continue across several blocks.

After $s$ repairs, an ordinary output prefix $B_n$ determines a list
of at most four words containing $A_{n+H_s}$. The ordinary code avoids
the markers used to announce repairs. Each completed repair stores
an exact source prefix, so only the current ordinary segment contributes
a list; ambiguities do not accumulate. Thus
$K(A_{n+H_s})\leq K(B_n)+O(1)$ at ordinary lengths, uniformly in the state.

The finite
\emph{reserve} $H_s$ means that the encoded source extends beyond
the current output length. We need this reserve because
Lemma~\ref{lem:sparse} allows a loss of $H_s/2$ in the complexity
gained during a repair. When only a few new output bits have been
written, the estimate may give no positive gain. Moreover, $H_s$
grows with $s$, so this loss cannot be absorbed into a fixed $O(1)$
term. The extra source bits cover the beginning of a repair and
compensate for this loss at later lengths.

A repair starting at $a$ writes a shortest program $p$ for
$A_{M+H_{s+1}}$ conditional on $(B_a,K(B_a))$. A short bridge before
$p$ carries a marker and the data needed to locate and run $p$,
including $K(B_a)$. The incompressibility of program prefixes,
together with a suitable choice of free bits, gives
\[
 K(B_{a+t})\geq K(B_a)+\left(1-\frac1{k_s}\right)t-\frac{H_s}{2}
 \qquad(0\leq t\leq M-a).
\]

To choose $a$, put
$C_j(x)=K(x\mid B_N,K(B_N))$, and choose $a$ to maximize
\[
 C_j(A_{m+H_s})-\left(1-\frac1{k_s}\right)(m-N)
 \qquad(N\leq m\leq M).
\]
If $M-a$ is sufficiently large, we perform a repair there; otherwise
we continue ordinary coding to $M$. For $H_s\leq t\leq M-a$,
ordinary decoding, the chain rule, and maximality at $m=a+t-H_s$ give
\[
 K(A_{a+t})\leq K(B_a)+\left(1-\frac1{k_s}\right)(t-H_s)+O(1).
\]
Since $k_s\geq2$, the saved term $(1-1/k_s)H_s$ is at least
$H_s/2$. Comparing with the repair growth estimate gives
$K(A_{a+t})\leq K(B_{a+t})+O(1)$ throughout this part of the repair.
The same maximum condition bounds $|p|$, ensuring that the bridge
and program fit in the available interval.

For ordinary lengths and the first $H_s$ positions of a repair,
we use a delayed scan. If the stored source word is long enough,
truncate it to length $n$. Otherwise scan only up to $n-H_s$.
For $n\leq a+H_s$, this is at most $a$, so the scan avoids the new
repair marker and ordinary decoding gives a list containing $A_n$.
Earlier repairs are decoded from the full input $B_n$. The list has
at most four entries, proving the comparison at these remaining
lengths and hence $A\leq_K B$.

To obtain dimension one, the maximum property also gives a lower
bound on the complexity gained in the ordinary part of each block.
Together with the repair estimate, it yields growth at rate
$1-1/k_s$, up to an error of order $(k_s+1)^2\log(M+2)$.
If repairs stopped, reverse bounded-list decoding for the final
ordinary segment would give $\dimH(B)\leq\dimH(A)$. The block
estimates would instead give
$\dimH(B)\geq1-1/k_s\geq1-1/k_*>\dimH(A)$.
Thus repairs occur infinitely often and $k_s\to\infty$. The slow
growth of $k_s$ keeps the accumulated errors sublinear, and summing
the block estimates gives $K(B_N)\geq N-o(N)$.

Finally, these infinitely many repairs give $A\leq_T B$.
From $B$, read successive repair records and run their programs with
condition $(B_a,K(B_a))$, using the recorded complexity value.
The programs return entire source prefixes of unbounded lengths,
which determine $A$.

\subsection{Construction}
For $s\in\N$, put
\[
 k_s=k_*+\log(s+2)-1,\qquad
 H_s=2c_*+4\log k_s,
\]
where $c_*$ is a positive integer to be fixed in Lemma~\ref{lem:sparse}.
The following estimate controls complexity throughout a bridge whose
prescribed positions lie in one residue class.

\begin{lemma}\label{lem:sparse}
For all sufficiently large $c_*$, the following holds for every
$s\in\N$ and $\sigma\in\Strings$.
Fix $r\in\{0,\ldots,4k_s-1\}$ and $z\in\Cantor$. Prescribe $Y(r+4k_sj)=z(j)$ for all $j$,
and choose the bits $Y(i)$ with $i\not\equiv r\pmod{4k_s}$ fairly and
independently. With probability at least
$3/4$,
\[
 K(\sigma Y_t)\geq K(\sigma)+\left(1-\frac1{k_s}\right)t-\frac{H_s}{2}, \quad\forall t\geq0.
\]
\end{lemma}

\begin{proof}
Call $\tau\in2^t$ \emph{compatible} if $\tau(r+4k_sj)=z(j)$
for every $j\in\N$ with $r+4k_sj<t$.
Among the first $t$ positions, at most $t/(4k_s)+1$ are prescribed.
Consequently for every compatible $\tau\in2^t$,
\begin{equation}\label{eq:trackprobability}
 \Pr(Y_t=\tau)\leq 2^{-(1-1/(4k_s))t+1}.
\end{equation}
Put
\[
 q(\sigma)=\frac12\sum_{t\geq0}\frac1{(t+1)^2}
                    \sum_{|\tau|=t}2^{-K(\sigma\tau)}.
\]
This is a left-r.e.\ discrete semimeasure: for each $t$, the words
$\sigma\tau$ with $|\tau|=t$ have unique decompositions, so
$\sum_\sigma q(\sigma)
 \leq\frac12\sum_{t\geq0}(t+1)^{-2}\leq1$.
By~\eqref{eq:semimeasure}, $q(\sigma)\leq c_0 2^{-K(\sigma)}$ for some
constant $c_0$ independent of $s$ and $\sigma$.
Putting $u=t/k_s$ and using $k_s\geq1$, we have
\[
 (t+1)^2 2^{-t/(4k_s)}
 \leq k_s^2(u+1)^2 2^{-u/4}=O(k_s^2).
\]
Thus $2^{-t/(4k_s)}\leq O(k_s^2)(t+1)^{-2}$ uniformly in $t$ and $s$.
Using the definition of $q$, we obtain
\begin{equation}\label{eq:weightedkraft}
 \begin{aligned}
 \sum_{t\geq0}2^{-t/(4k_s)}
       \sum_{|\tau|=t}2^{-K(\sigma\tau)}
 &\leq O(k_s^2)\sum_{t\geq0}\frac1{(t+1)^2}
       \sum_{|\tau|=t}2^{-K(\sigma\tau)}\\
 &=O(k_s^2)q(\sigma)
 \leq c_1 k_s^2 2^{-K(\sigma)},
 \end{aligned}
\end{equation}
where $c_1$ is independent of $s$ and $\sigma$.

Call $\tau\in2^t$ \emph{bad} if
$K(\sigma\tau)<K(\sigma)+(1-1/k_s)t-\frac{H_s}{2}$.
For each compatible bad word, \eqref{eq:trackprobability} and the
definition of badness give
\[
 \begin{aligned}
 \Pr(Y_t=\tau)
 &\leq 2^{1-(1-1/(4k_s))t}\\
 &\leq 2^{1-(1-1/k_s)t-t/(4k_s)}\\
 &<2^{K(\sigma)+1-\frac{H_s}{2}}
       2^{-t/(4k_s)}2^{-K(\sigma\tau)},
 \end{aligned}
\]
where the last inequality uses badness.
Summing over all bad words and lengths, and
using~\eqref{eq:weightedkraft}, bounds the probability of any bad
prefix by
\[
 c_1k_s^2 2^{1-\frac{H_s}{2}}
 =c_1k_s^2 2^{1-c_*-2\log k_s}
 \leq c_1 2^{1-c_*}\leq\frac14,
\]
where the last inequality holds when $c_*$ is sufficiently large.
\end{proof}

A \emph{state} is specified by an endpoint $E$, a repair count $s$, and the stored word $A_{E+H_s}$. We denote this state by $(E,s)$, with its stored word understood. Its set of control positions is
\[
 \mathcal C_{E,s}=\{E+4k_sr:r\in\N\}.
\]
The initial state is $(E,s)=(0,0)$, with stored word $A_{H_0}$.
For a tail $Y$ placed at position $E$ and $u\in\mathcal C_{E,s}$,
the marker event at $u$ is
\[
 \begin{split}
 \mathcal M_{E,s,u}=\{Y\in\Cantor:
 &Y(u-E)=Y(u-E+4k_s(2\log(u+2)+1))=1,\\
 &Y(u-E+4k_si)=0\quad(1\leq i\leq 2\log(u+2))\}.
 \end{split}
\]
Thus the pattern on consecutive control positions is $1\,0^{2\log(u+2)}\,1$.
Only starts in $\mathcal C_{E,s}$ are considered. Define the marker-free
closed class and its conditional measure by
\[
 P_{E,s}=\Cantor\setminus\bigcup_{u\in\mathcal C_{E,s}}\mathcal M_{E,s,u},
 \qquad
 \nu_{E,s}(Q)=\frac{\lambda(Q\cap P_{E,s})}{\lambda(P_{E,s})}
\]
for Borel $Q\subseteq\Cantor$.
\begin{lemma}\label{lem:markerprob}
For every $E,s\in\N$, $\nu_{E,s}$ is a recursive probability measure,
uniformly in $E,s$, with support contained in $P_{E,s}$ and $\nu_{E,s}[\tau]\leq2\cdot2^{-|\tau|}$ for all $\tau$.
\end{lemma}
\begin{proof}
A marker starting at $u$ has fair-coin probability $2^{-2\log(u+2)-2}$.
Thus
\[
 \lambda(P_{E,s})\geq1-\sum_{u\geq E}2^{-2\log(u+2)-2}
 \geq1-\frac14\sum_{u\geq E}(u+2)^{-2}>\frac34.
\]
The marker events with starts below a given cutoff form a uniformly
recursive clopen set, and the remaining events have measure bounded
by an effectively convergent tail of this series. Hence
$\lambda(P_{E,s}\cap[\tau])$ is recursive uniformly in $E,s,\tau$. Dividing by
$\lambda(P_{E,s})>1/2$ gives the asserted recursiveness and upper
bound. Since $P_{E,s}$ is closed and has full $\nu_{E,s}$-measure,
it contains the support.
\end{proof}

In state $(E,s)$, the word $A_{E+H_s}$ is already stored, so the
remaining source sequence to be encoded is $A[E+H_s:\infty)$. Put
$X^{E,s}=A[E+H_s:\infty)$, and let $Y^{E,s}$ be obtained by the interval
construction in the proof of Lemma~\ref{lem:listcode}, applied to
$X^{E,s}$ and $\nu_{E,s}$.

\begin{lemma}\label{lem:markermeasure}
For every $E,s\in\N$, $\nu_{E,s}[\tau]>0$ exactly when $\tau$ contains
no completed marker for the state $(E,s)$, and then
$\nu_{E,s}[\tau]\geq\frac14\,2^{-|\tau|}$.
For each fixed state $(E,s)$, the source sequence
$X^{E,s}=A[E+H_s:\infty)$ and the corresponding output sequence
$Y^{E,s}$ chosen above satisfy $K(Y^{E,s}_t)\leq K(X^{E,s}_t)+O(1)$,
where the constant depends only on $E,s$.
\end{lemma}
\begin{proof}
Suppose that $\tau\in2^t$ contains no completed marker, and put
$b=E+t$. Among markers starting before $b$, at most one remains
compatible with $\tau$: its initial $1$ must be the last $1$ on the
control positions present in $\tau$, followed only by zeros. Its
conditional probability given $[\tau]$ is at most $1/2$, since at
least one prescribed bit lies beyond $\tau$. Markers starting at or
after $b$ depend only on the remaining bits, so
\[
 \lambda(P_{E,s}\mid[\tau])
 \geq1-\frac12-\sum_{u\geq b}2^{-2\log(u+2)-2}
 >\frac14.
\]
Consequently
$\nu_{E,s}[\tau]\geq\lambda(P_{E,s}\cap[\tau])
\geq2^{-t}/4$. A prefix containing a completed marker has $\nu_{E,s}$-measure
zero. This also makes positivity of $\nu_{E,s}[\tau]$ uniformly
decidable.

For the complexity bound, fix $E,s$ and use the intervals $J^{E,s}_\tau$
and their rational enlargements $\widetilde J^{E,s}_\tau$ from the proof of
Lemma~\ref{lem:listcode}, with $\varepsilon=2^{-t-4}$. On input
$\xi\in2^t$, compute the list
\[
 \{\tau\in2^t:\nu_{E,s}[\tau]>0\text{ and }
       \widetilde J^{E,s}_\tau\cap I_\xi\neq\varnothing\}.
\]
This list is recursive. Each selected $J^{E,s}_\tau$ meets
$I_\xi+[-2\varepsilon,2\varepsilon]$, an interval of length
$5\cdot2^{-t}/4$. The positive intervals $J^{E,s}_\tau$ have disjoint
interiors and length at least $2^{-t}/4$, so at most seven can meet
this interval. By the choice of $Y^{E,s}$, $x_{X^{E,s}}\in I_{X^{E,s}_t}\cap J^{E,s}_{Y^{E,s}_t}$ and $\nu_{E,s}[Y^{E,s}_t]>0$.
Thus the list computed from $X^{E,s}_t$ contains $Y^{E,s}_t$, and
\eqref{eq:boundedlist} gives the desired complexity bound.
\end{proof}

Suppose the current state is $(E,s)$ and ordinary coding has produced the
prefix $B_a$. We describe a repair that starts at position $a$ and ends at
the next block endpoint $M>a$, extending $B_a$ to $B_M$ by filling the interval $[a,M)$.
Write
\[
 \sigma=B_a,\qquad T=M-a,\qquad
 v=A_{a+H_s},\qquad w=A_{M+H_{s+1}}.
\]
For a positive integer $\kappa_1$, whose value will be fixed before the stage instructions, set
\[
 h'=\kappa_1\log(M+2),\qquad L=4k_s h'.
\]
The repair has the following format.

\smallskip
\noindent\emph{The bridge.}
The interval $[a,a+L)$ is the bridge. Let $u_0=E+4k_s\left\lceil\frac{a-E}{4k_s}\right\rceil \in \mathcal C_{E,s}$. The control positions are
$u_r=u_0+4k_sr$ ($0\leq r<h'$).
Prescribe $B(u_0)=1$. If this completes a marker starting at some $u$ with $E\leq u<a$,
use that completion. Otherwise prescribe 
\[
 B(u_0+4k_si)=0\quad(1\leq i\leq 2\log(u_0+2)),\qquad
 B(u_0+4k_s(2\log(u_0+2)+1))=1.
\] 
In the latter case the $1$ at $u_0$ interrupts every older unfinished run
whose start is relevant to the current state, so this prescription produces the first marker completion after $a$.

Immediately after the first marker completion, write $a$, $K(\sigma)$,
$L$, and $|p|$, in that order, on successive control positions,
using a fixed self-delimiting integer code of logarithmic length.
These four integers form the \emph{metadata}.
Fill every remaining control position in the bridge with $0$.
All bridge positions outside $\mathcal C_{E,s}$ are free.
Lemma~\ref{lem:feasible} below ensures that the marker and the four
integer codes fit in the bridge.

\begin{figure}[htbp]
\centering
\begin{tikzpicture}[x=1cm,y=1cm,font=\footnotesize,
                    line width=0.45pt,>=latex]
  \draw[<->] (0,1.57) -- node[above] {Repair length $T=M-a$} (11.8,1.57);
  \draw[fill=black!7] (0,0) rectangle (4.8,1.02);
  \draw[fill=black!20] (4.8,0) rectangle (7.8,1.02);
  \draw[fill=white] (7.8,0) rectangle (11.8,1.02);
  \node at (2.4,0.68) {Bridge};
  \node at (2.4,0.29) {$L=4k_s h'$ bits};
  \node at (6.3,0.68) {Program $p$};
  \node at (6.3,0.29) {$|p|$ bits};
  \node at (9.8,0.68) {Free bits};
  \node at (9.8,0.29) {$T-L-|p|$ bits};
  \foreach \x in {0,4.8,7.8,11.8}
    \draw (\x,0) -- (\x,-0.13);
  \node[anchor=north west,inner sep=1pt] at (0,-0.18) {$a$};
  \node[anchor=north,inner sep=1pt] at (4.8,-0.18) {$a+L$};
  \node[anchor=north,inner sep=1pt] at (7.8,-0.18) {$a+L+|p|$};
  \node[anchor=north east,inner sep=1pt] at (11.8,-0.18) {$M$};

  \node[anchor=west,inner sep=0pt] at (0,-0.94)
    {Bridge detail: $h'$ prescribed positions and $L-h'$ free positions};
  \draw[fill=white] (0,-1.75) rectangle (11.8,-1.32);
  \foreach \x in {0.40,2.45,4.50,8.60,10.65}
    \draw[fill=black!75] (\x,-1.75) rectangle +(0.25,0.43);
  \node[fill=white,inner sep=3pt] at (6.675,-1.535) {$\cdots$};
  \draw (0.525,-1.75) -- (0.525,-1.88);
  \draw (2.575,-1.75) -- (2.575,-1.88);
  \draw (10.775,-1.75) -- (10.775,-1.88);
  \node[anchor=north,inner sep=1pt] at (0.525,-1.92) {$u_0$};
  \node[anchor=north,inner sep=1pt] at (2.575,-1.92) {$u_0+4k_s$};
  \node[anchor=north east,inner sep=1pt] at (11.8,-1.92)
    {$u_0+4k_s(h'-1)$};
  \draw[fill=black!75] (3.0,-2.66) rectangle +(0.20,0.20);
  \node[anchor=west,inner sep=0pt] at (3.36,-2.56) {prescribed};
  \draw[fill=white] (6.9,-2.66) rectangle +(0.20,0.20);
  \node[anchor=west,inner sep=0pt] at (7.26,-2.56) {free};
\end{tikzpicture}
\caption{Position allocation within one repair (not to scale). The prescribed
bridge bits complete the marker, encode the metadata, and provide any zero
padding.}
\label{fig:repair-layout}
\end{figure}
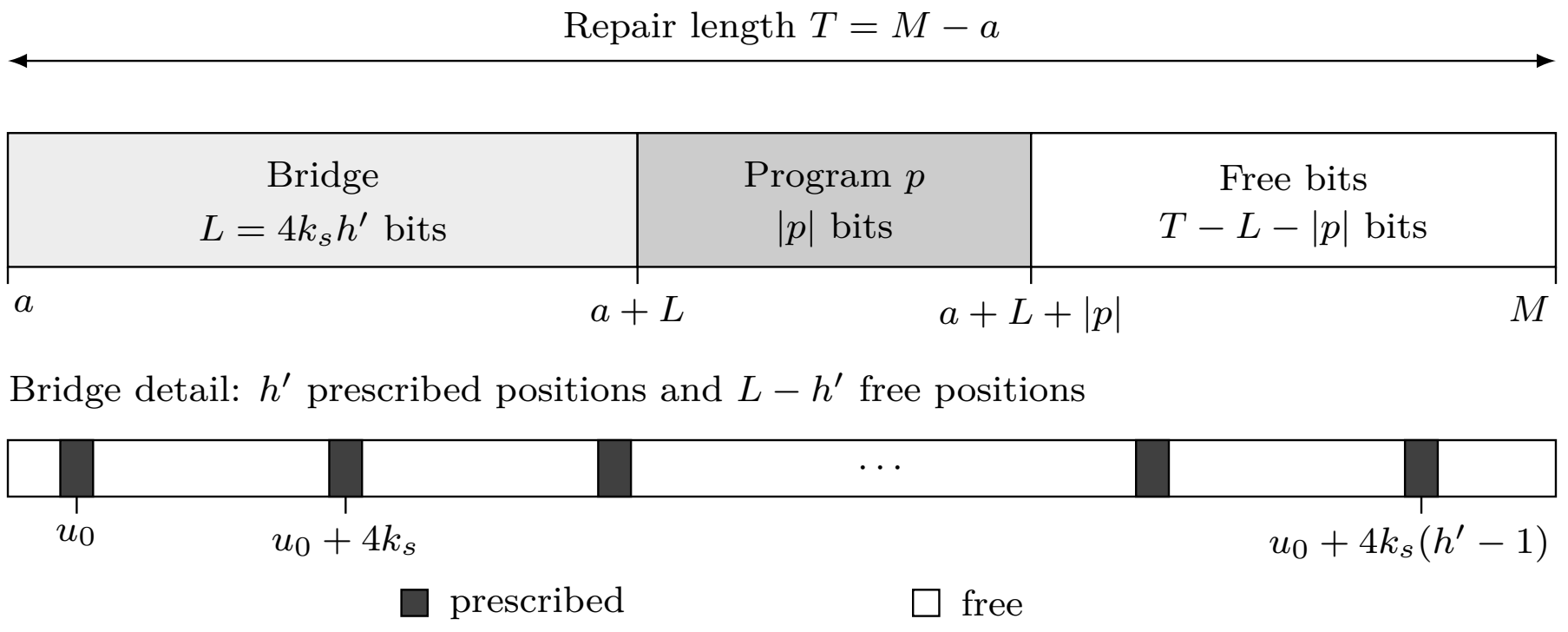

\smallskip
\noindent\emph{The program and the final free positions.} Let $p$ be the canonical shortest program for $w$ conditional on $(\sigma,K(\sigma))$.
Write $B(a+L+i)=p(i)$ for $0\leq i<|p|$.
The full set of free positions is $\bigl([a,a+L)\setminus\mathcal C_{E,s}\bigr)
 \cup[a+L+|p|,M)$.
Subject to the feasibility inequalities proved below, the three parts
satisfy $a<a+L\leq a+L+|p|\leq M$.
Figure~\ref{fig:repair-layout} shows how these positions are allocated.

Choose a sufficiently large positive integer $\kappa_1$ as specified in
Lemma~\ref{lem:feasible}.

We can now give the \textbf{stage instructions}.
At stage $j$, set $N=2j^2$ and $M=2(j+1)^2$. $B_N$ is already defined, in a state $(E,s)$ with $E\leq N$.
Define
\[
 \Theta_j=4\kappa_1(k_s+1)^2\log(M+2).
\]
Put
\[
 C_j(x)=K(x\mid B_N,K(B_N)),
\]
and set
\[
 f_j(m)=C_j(A_{m+H_s})-\left(1-\frac1{k_s}\right)(m-N).
\]
Define
\[
 a=\min\{m\in[N,M]:
               f_j(m)=\max_{u\in[N,M]}f_j(u)\},\qquad T=M-a.
\]
In every block, the maximum property gives, for $a\leq m\leq M$,
\begin{align}
 C_j(A_{m+H_s})-C_j(A_{a+H_s})\le \left(1-\frac1{k_s}\right)(m-a),
 \label{eq:maxright}\\
 C_j(A_{a+H_s})\ge C_j(A_{N+H_s})+\left(1-\frac1{k_s}\right)(a-N)
 \geq \left(1-\frac1{k_s}\right)(a-N).\label{eq:maxleft}
\end{align}

If $T<\Theta_j$, set $B_M=B_N\,Y^{E,s}[N-E:M-E)$
and retain state $(E,s)$.

If $T\geq\Theta_j$, first set $B_a=B_N\,Y^{E,s}[N-E:a-E)$.
Use the repair format above on $[a,M)$. Retaining all prescribed
marker, metadata, and program bits, choose the lexicographically least
extension $B_M\in2^M$ of $B_a$ satisfying
\begin{equation}\label{eq:repairgrowth}
 K(B_{a+t})\geq K(B_a)+\left(1-\frac1{k_s}\right)t-\frac{H_s}{2}
 \qquad(0\leq t\leq T).
\end{equation}
Enter state $(M,s+1)$, with stored word $A_{M+H_{s+1}}$.

This completes the stage instructions. We will prove that every proposed
repair permits a choice of free bits satisfying~\eqref{eq:repairgrowth}. The output is the union $B=\bigcup_{j\in\N}B_{2j^2}$.

\subsection{Verification}\label{sec:verification}
We prove that the construction is well defined and that its output $B$
satisfies $A\leq_K B$, $A\leq_T B$, and $\dimH(B)=1$.
We first justify decoding and
show that every proposed repair can be carried out. We then establish
$K(A_n)\leq K(B_n)+O(1)$ at every length. Finally, we estimate the
complexity gained in each block, show that repairs occur infinitely often,
and sum the block estimates to obtain $\dimH(B)=1$.
The infinitely many repair records also give $A\leq_T B$.
Unless explicitly stated otherwise, all constants below are uniform in the
block and the finite history.

A \emph{finite valid history} consists of finitely many completed blocks
following the stage instructions from the initial state, possibly followed
by an initial part of the ordinary coding in the next block. In every
completed repair, the marker and metadata fit in the bridge, the program
fits in the repair interval, all prescribed bits are written correctly,
and the free bits satisfy~\eqref{eq:repairgrowth}. The states and stored
words are updated as prescribed. The initial state alone is also a finite
valid history.

\begin{lemma}\label{lem:ordinary}
Let $B_n$ be a prefix of a finite valid history,
and suppose that $n$ lies outside the open interval $(a,M)$ for every
repair interval $[a,M)$ in that history. Let $s$ be the number of repairs
completed by length $n$. Thus $n$ may be any length during ordinary
coding, including a repair starting point, or the endpoint of a completed
repair; at a repair endpoint, that repair is included in $s$.

There is one partial recursive procedure which, given $B_n$ alone,
outputs at most 4 words of length $n+H_s$, one of which is
$A_{n+H_s}$. Consequently there is a constant $c_f$, independent of
$n$, $s$, and the history, such that $K(A_{n+H_s})\leq K(B_n)+c_f$
at every such length $n$.
\end{lemma}
\begin{proof}
The procedure reads $n$ as the length of its input. Start with state
$(0,0)$ and the stored word $A_{H_0}$. In the current state, scan the
control positions present in the input for the first marker completion.
Ordinary coding contains no marker. Since $n$ is not strictly inside a
repair, any marker found belongs to a repair whose endpoint is at most
$n$, so its entire record and program are present. Read the metadata,
recover $B_a$ by truncation, and run the recorded program with condition
$(B_a,K(B_a))$, using the value of $K(B_a)$ supplied by the metadata.
The program returns the exact new stored word. Recover the next block
endpoint from $a$ as
$M=2\left(\left\lfloor\sqrt{a/2}\right\rfloor+1\right)^2$.
Jump to $M$, increase the repair count, and repeat.
When no marker completion remains, the recovered state is $(E,s)$,
where $E$ is the last completed repair endpoint. The stored word is $A_{E+H_s}$, and the remaining suffix
satisfies $B[E:n)=Y^{E,s}\restr(n-E)$.
If $n=E$, return the stored word as a singleton list. Otherwise apply
the decoder of Lemma~\ref{lem:listcode} to $B[E:n)$, using the uniformly recursive
index for $\nu_{E,s}$ supplied by Lemma~\ref{lem:markerprob}.
Prepend $A_{E+H_s}$ to each candidate. This gives at most four words
of length $n+H_s$, one of which is $A_{n+H_s}$.
The same procedure works for all the stated inputs, so
\eqref{eq:boundedlist} gives the uniform constant $c_f$.
\end{proof}

\begin{lemma}\label{lem:freebits}
Fix a finite valid history ending with an ordinary prefix $\sigma=B_a$
in state $(E,s)$, and the data of a proposed repair ending at the next
block endpoint $M$. Suppose that $L+|p|\leq T=M-a$ and that the marker and metadata fit in the bridge. Then the free bits can be chosen so that~\eqref{eq:repairgrowth} holds.
\end{lemma}
\begin{proof}
Choose all free bits independently with probability $1/2$ for each value. In the bridge, only the bits on one residue class modulo $4k_s$ are prescribed. Their values, including the metadata, were fixed before this choice. Lemma~\ref{lem:sparse} implies~\eqref{eq:repairgrowth} for every $0\leq t\leq L$ with probability at least $3/4$.

For $L\leq t\leq T$, set
\[
 j(t)=\min\{t-L,|p|\},\qquad r(t)=t-h'-j(t).
\]
Let $\eta$ be the word formed by the free bits of this repair, in their order of occurrence. Thus $\eta_{r(t)}$ consists of the free bits already read at length $a+t$. From $B_{a+t}$, a fixed partial recursive procedure recovers $(\sigma,p_{j(t)},\eta_{r(t)})$.
As in the proof of Lemma~\ref{lem:ordinary}, scan successive marker
completions and read their metadata, passing to the next state only when
the repair endpoint is strictly smaller than the input length.
Stop at the first repair whose endpoint is at least the input length;
its entire bridge is present since $t\geq L$.
The recovered state and metadata determine $a$, $L$, $|p|$, and the
free positions, so truncation gives $\sigma$, and the program bits and
free bits already present give $p_{j(t)}$ and $\eta_{r(t)}$, respectively.

For a fixed $t$, the pair $(\sigma,p_{j(t)})$ is fixed, and $\eta_{r(t)}$ is uniform on $2^{r(t)}$. Conditional prefix-free counting gives
\[
 \Prob\bigl(K(\eta_{r(t)}\mid(\sigma,p_{j(t)})^*)<r(t)-z\bigr)\leq2^{-z}.
\]
Choose $z=\log(T+1)+2$. By the union bound,
\[
 \begin{aligned}
 &\Prob\bigl(\exists t\in[L,T]\cap\N:\quad
 K(\eta_{r(t)}\mid(\sigma,p_{j(t)})^*)<r(t)-z\bigr)\\
 &\qquad\leq\sum_{t=L}^{T}2^{-z}
 \leq(T+1)2^{-\log(T+1)-2}\leq\frac14,
 \end{aligned}
\]
since $\log(T+1)=\lceil\log_2(T+1)\rceil$.
Hence, with probability at least $3/4$, simultaneously for $L\leq t\leq T$,
\[
 K(\eta_{r(t)}\mid(\sigma,p_{j(t)})^*)
 \geq r(t)-\log(T+1)-2.
\]
Fix $0\leq i\leq|p|$, and let $q$ be a shortest description of $p_i$
conditional on $(\sigma,K(\sigma))$. With the same condition, a fixed
machine first reads $q$, recovers $p_i$, and then runs the original
machine on $p_i$ followed by the remaining input. Its domain is
prefix-free, since $q$ is self-delimiting and the reconstructed program
must belong to the original prefix-free domain. Thus $q\,p[i:|p|)$
describes $w$, and the minimality of $p$ gives
\[
 |p|\leq |q|+|p|-i+O(1).
\]
Since $|q|=K(p_i\mid\sigma,K(\sigma))$, rearranging yields
$K(p_i\mid\sigma,K(\sigma))\geq i-O(1)$, with a constant independent
of $i$, $p$, and $\sigma$.
Thus, for every such assignment, the chain rule gives
\begin{align*}
 K(B_{a+t})
 &\geq K(\sigma,p_{j(t)},\eta_{r(t)})-O(1)\\
 &\geq K(\sigma,p_{j(t)})+
       K(\eta_{r(t)}\mid(\sigma,p_{j(t)})^*)-O(1)\\
 &\geq K(\sigma)+j(t)+r(t)-\log(T+1)-\kappa_2\\
 &=K(\sigma)+t-h'-\log(T+1)-\kappa_2.
\end{align*}
Here $\kappa_2$ is independent of $\kappa_1$. Our fixed choice
$3\kappa_1\geq1+\kappa_2$, together with $h'=\kappa_1\log(M+2)$, gives $4h'\geq h'+\log(T+1)+\kappa_2$.
For $t\geq L=4k_s h'$, we have $t/k_s\geq4h'$, so the last lower bound is at least $K(\sigma)+(1-1/k_s)t$. Thus~\eqref{eq:repairgrowth} holds throughout $L\leq t\leq T$ with probability at least $3/4$. The two events have intersection of probability at least $1/2$, so the free bits can be chosen to satisfy~\eqref{eq:repairgrowth}.
\end{proof}

\begin{lemma}\label{lem:feasible}
Let $\kappa_2$ be the constant defined in the proof of
Lemma~\ref{lem:freebits}, and let $\kappa_0\geq0$ be a sufficiently large
integer. Then
fix any sufficiently large positive integer $\kappa_1$ satisfying
$3\kappa_1\geq1+\kappa_2$ and $4\kappa_1\geq2c_*+4+\kappa_0$.
Fix a finite valid history ending at $N=2j^2$, and put
$M=2(j+1)^2$. If the stage instructions propose a repair on $[a,M)$
in block $j$, put
\[
 \sigma=B_a,\qquad v=A_{a+H_s},\qquad w=A_{M+H_{s+1}}.
\]
Then
\begin{align}
 K(B_N)+C_j(v)&\leq K(B_a)+O(1),\label{eq:sourcepair}\\
 |p|&\leq C_j(w)-C_j(v)+\kappa_0\log(M+2),
 \label{eq:programbound}\\
 L+|p|+H_s&\leq M-a.\label{eq:fits}
\end{align}
These estimates are uniform in the block and the finite valid history.
Moreover, every finite valid history ending at a block endpoint can be
extended through the next block to a finite valid history. Consequently,
every stage of the construction is well defined.
\end{lemma}
\begin{proof}
Fix such a finite valid history, and let $(E,s)$ be its current state.
Suppose a repair is proposed on $[a,M)$. Extend the ordinary code to $B_a$,
still obtaining a finite valid history, and put
$T=M-a$.
The repair threshold gives $M\geq T\geq4\kappa_1(k_s+1)^2\log(M+2)$, so
$k_s+1\leq M$. Since $k_{s+1}\leq k_s+1$, we have
\[
 \begin{aligned}
 H_s&\leq H_{s+1}\leq2c_*+4\log(k_s+1)\\
    &\leq(2c_*+4)\log(M+2)\leq4\kappa_1\log(M+2)\leq T,
 \end{aligned}
\]
so $|v|,|w|\leq2M$.

By Lemma~\ref{lem:ordinary}, applied to this history ending at $a$,
$\sigma$ determines a list of at most four words containing $v$. Since $N\le a<M$, $N=2\left\lfloor\sqrt{a/2}\right\rfloor^2$
is recursive from $|\sigma|$. Thus $(B_N,v)$ belongs to a fixed
bounded list obtained from $\sigma$, and the chain rule proves
\eqref{eq:sourcepair}.

Since
$N\leq M$ and $|v|\leq2M$, we have $K(B_N)\leq N+O(\log(N+2))=O(M)$ and 
\[
 \begin{aligned}
 C_j(v)=K(v\mid B_N,K(B_N))\leq K(v)+O(1)
       =O(M).
 \end{aligned}
\]
Each nonnegative integer $u$ has a self-delimiting code of length
$O(\log(u+2))$. Thus the integers $K(B_N)$ and $C_j(v)$ can each
be specified using $O(\log(M+2))$ bits. Together with the index of
$v$, these codes have total length $O(\log(M+2))$ and allow us to
recover $(B_N,K(B_N),v,C_j(v))$ from $(\sigma,K(\sigma))$.
 Since $p$ is a shortest
description of $w$ conditional on $(\sigma,K(\sigma))$, we obtain
\[
 \begin{aligned}
 |p|=K(w\mid\sigma,K(\sigma))\leq K(w\mid B_N,K(B_N),v,C_j(v))+O(\log(M+2)).
 \end{aligned}
\]
Apply the conditional chain rule~\eqref{eq:chain} with $x=v$,
$y=w$, and $z=(B_N,K(B_N))$. Since
$K(v\mid B_N,K(B_N))=C_j(v)$, this gives
\[
 \begin{aligned}
 K(v,w\mid B_N,K(B_N))
 =C_j(v)
+K(w\mid B_N,K(B_N),v,C_j(v))+O(1).
 \end{aligned}
\]
Finally, $v=A_{a+H_s}$ is a prefix of $w=A_{M+H_{s+1}}$.
Prepend a self-delimiting code for $|v|$ to a shortest description
of $w$ conditional on $(B_N,K(B_N))$. After recovering $w$,
truncate it to length $|v|$ to recover $v$, and output $(v,w)$.
Since $|v|\leq2M$, the additional integer code has length
$O(\log(M+2))$, so
\[
 \begin{aligned}
 K(v,w\mid B_N,K(B_N))
 &\leq K(w\mid B_N,K(B_N))+O(\log(M+2))\\
 &=C_j(w)+O(\log(M+2)).
 \end{aligned}
\]
Combining these estimates yields
\[
 \begin{aligned}
 |p|&\leq K(v,w\mid B_N,K(B_N))-C_j(v)+O(\log(M+2))\\
    &\leq C_j(w)-C_j(v)+O(\log(M+2)).
 \end{aligned}
\]
The constants depend only on the fixed decoders and complexity
conventions, so~\eqref{eq:programbound} holds for sufficiently large
$\kappa_0$.

Since $a$ maximizes $f_j$, we have 
\[
 C_j(A_{M+H_s})-C_j(v)
 \leq\left(1-\frac1{k_s}\right)T.
\]
Also, $w=A_{M+H_{s+1}}$ extends $A_{M+H_s}$ by at most four bits,
since $H_{s+1}-H_s\leq4$. The length and contents of this suffix
can be specified using $O(1)$ bits. Appending this information to
a description of $A_{M+H_s}$, with the same condition
$(B_N,K(B_N))$, gives $C_j(w)\leq C_j(A_{M+H_s})+O(1)$.
This constant error is absorbed into
the preceding $O(\log(M+2))$ error. Choose $\kappa_0$ sufficiently large
to bound the error terms in both~\eqref{eq:programbound} and the
following estimate. Then
\begin{equation}\label{eq:programrate}
 \begin{aligned}
 |p|&\leq C_j(w)-C_j(v)+O(\log(M+2))\\
    &\leq C_j(A_{M+H_s})-C_j(v)+O(\log(M+2))\\
    &\leq\left(1-\frac1{k_s}\right)T+\kappa_0\log(M+2).
 \end{aligned}
\end{equation}
The repair threshold and our choice of $\kappa_1$ give
\[
 \begin{aligned}
 \frac{T}{k_s}&\geq4\kappa_1(k_s+1)\log(M+2)=L+4\kappa_1\log(M+2)\\
              &\geq L+(2c_*+4+\kappa_0)\log(M+2)
              \geq L+H_s+\kappa_0\log(M+2).
 \end{aligned}
\]
Together with~\eqref{eq:programrate}, this proves~\eqref{eq:fits}.

It remains to fit the marker and metadata in the bridge. The
integers $a,L,|p|$ are at most $M$ by~\eqref{eq:fits}, and
$K(\sigma)\leq a+O(\log(a+2))$. The first control position after $a$
is less than $a+4k_s\leq a+L \le M$ by \eqref{eq:fits}. The marker and the four self-delimiting
integer codes therefore occupy $O(\log(M+2))$ control positions.
Since the implicit constant is independent of $\kappa_1$ and of
the block and history, these encodings fit in the
$h'=\kappa_1\log(M+2)$ available positions for every sufficiently large
$\kappa_1$.

Lemma~\ref{lem:freebits} supplies a choice of free bits satisfying
\eqref{eq:repairgrowth}, so the lexicographically least such choice exists.
The recorded program supplies the next stored word. Thus the history
extends through block $j$ to a finite valid history. Starting from the
initial state, induction now shows that every stage can be carried out.
\end{proof}

Lemma~\ref{lem:feasible} shows that the infinite output $B$ exists.
For the remainder of this section, fix this $B$. All states and repair
intervals below refer to its construction. We first prove the complexity
comparison at every length, including the transition from ordinary
coding to a repair.
\begin{lemma}\label{lem:delayed}
There is a constant $c_\ell$ such that
$K(A_n)\leq K(B_n)+c_\ell$
at every length $n$ satisfying at least one of the following conditions:
\begin{enumerate}
\item $n\notin(a,M)$ for every repair interval $[a,M)$;
\item $a\leq n\leq a+H_s$ for some repair interval $[a,M)$ starting
in state $(E,s)$.
\end{enumerate}
In fact, there is one partial recursive procedure which, given any
such prefix $B_n$ alone, outputs a list of at most four words of
length $n$ containing $A_n$.
\end{lemma}
\begin{proof}
Given a word of length $n$, begin in the initial state. In a decoded state $(E,s)$, with stored word $A_{E+H_s}$, act as follows. If $n\leq E+H_s$, output the stored word truncated to length $n$. Otherwise set $q=n-H_s>E$ and scan only marker completions visible in the prefix of length $q$. If none is found, apply the ordinary inverse code to $B[E:q)$ and prepend the stored word. The resulting candidates have length $(E+H_s)+(q-E)=n$.
If a marker is found, use the full input of length $n$ to read and execute its repair record. Enter the new state and repeat. 

We verify that this procedure works at the stated lengths.
Suppose the input ends inside a repair interval $[a,M)$ starting in state $(E,s)$.
The hypothesis of the lemma gives $n\leq a+H_s$, hence $q=n-H_s\leq a$.
If the stored word is not already long enough, the extra $H_s$ stored source bits allow us to obtain the required list using only the ordinary code $B[E:q)$.
The delayed scan therefore cannot see this repair's marker, so every marker actually processed belongs to a completed repair whose full record is present in $B_n$.
Every earlier repair occupies an interval $[a',M')$ with $M'\leq n$.
If its bridge has length $L'$ and its old state has reserve $H_{s'}$,
then its marker completion length satisfies $c'\leq a'+L'$.
Thus~\eqref{eq:fits} gives
$n-H_{s'}\geq M'-H_{s'}\geq a'+L'\geq c'$.
Whenever the algorithm needs to scan past that old state, the marker is therefore visible, and the full record is present in its input. It can process the earlier repairs correctly. The same argument applies to an ordinary current prefix. As before, only the last inverse code produces a list, whose size is at most four. Equation~\eqref{eq:boundedlist} proves the required inequality.
\end{proof}

\begin{lemma}\label{lem:domination}
$A\leq_K B$.
\end{lemma}
\begin{proof}
Lemma~\ref{lem:delayed} covers ordinary positions, completed endpoints,
and the first $H_s$ positions of each repair. Consider the remaining
case $n=a+t$ with $H_s\leq t\leq T$, in block $j$.
Put $m=n-H_s$, so $a\leq m\leq M$.
The chain rule, \eqref{eq:maxright}, and~\eqref{eq:sourcepair} give
\[
 \begin{split}
 K(A_n)
 &\leq K(B_N)+C_j(A_n)+O(1)\\
 &\leq K(B_N)+C_j(A_{a+H_s})+\left(1-\frac1{k_s}\right)(t-H_s)+O(1)\\
 &\leq K(B_a)+\left(1-\frac1{k_s}\right)(t-H_s)+O(1)\\
 &\leq K(B_{a+t})+O(1).
 \end{split}
\]
The last step uses~\eqref{eq:repairgrowth} and $k_s\geq2$.
All constants are uniform.
\end{proof}

It remains to show that $B$ has dimension one. For each block
$[2j^2,2(j+1)^2]$, let $s_j$ be the number of repairs completed
before that block.

\begin{lemma}\label{lem:blockestimate}
There is a fixed constant $\kappa_E$ such that, on putting
\[
 E_s(M)=\kappa_E(k_s+1)^2\log(M+2),
\]
every block $[N,M]=[2j^2,2(j+1)^2]$, with $s=s_j$, satisfies
\begin{equation}\label{eq:blockestimate}
 K(B_M)-K(B_N)
 \geq\left(1-\frac1{k_s}\right)(M-N)-E_s(M).
\end{equation}
\end{lemma}

\begin{proof}
Write $\delta=M-N$ and $v=A_{a+H_s}$.
By~\eqref{eq:maxleft}, $C_j(v)\geq \left(1-\frac1{k_s}\right)(a-N)$. If a repair is made, \eqref{eq:sourcepair} and~\eqref{eq:repairgrowth} give
\begin{align*}
 K(B_M)
 &\geq K(B_a)+\left(1-\frac1{k_s}\right)(M-a)-\frac{H_s}{2}\\
 &\geq K(B_N)+C_j(v)+\left(1-\frac1{k_s}\right)(M-a)-\frac{H_s}{2}-O(1)\\
 &\geq K(B_N)+\left(1-\frac1{k_s}\right)\delta-\frac{H_s}{2}-O(1).
\end{align*}

If no repair is made, then $M-a<\Theta_j$.
Lemma~\ref{lem:ordinary} gives a bounded list containing
$A_{M+H_s}$ from $B_M$. Specifying $a\leq M$ allows each candidate
to be truncated to length $a+H_s$, while $N$ is recursively determined
by the endpoint $M$. Hence
$K(B_N,v)\leq K(B_M)+O(\log(M+2))$.
By the chain rule and \eqref{eq:maxleft},
\begin{align*}
 K(B_M)-K(B_N)
 &\geq C_j(v)-O(\log(M+2))\\
 &\geq \left(1-\frac1{k_s}\right)\delta-\Theta_j-O(\log(M+2)).
\end{align*}
Both cases yield \eqref{eq:blockestimate} after increasing
$\kappa_E$.
\end{proof}

\begin{lemma}\label{lem:infinite}
There are infinitely many repairs.
\end{lemma}

\begin{proof}
Suppose otherwise, and let $(E,s)$ be the final state. Then $B[E:n)=Y^{E,s}\restr(n-E)$ for all $n\geq E$.
By the final assertion of Lemma~\ref{lem:markermeasure}, applied to
these source and output tails, and with $E$, $s$, and $B_E$ fixed, $K(B_n)\leq K(A_{n+H_s})+O(1)$,
so $\dimH(B)\leq\dimH(A)$.

On the other hand, fix $j_0$ such that $s_j=s$ for all $j\geq j_0$.
Since $s$ is fixed, the factor $\kappa_E(k_s+1)^2$ is constant. For $J>j_0$,
\[
 \begin{aligned}
 \sum_{j=j_0}^{J-1}E_s(2(j+1)^2)
 &=\kappa_E(k_s+1)^2
   \sum_{j=j_0}^{J-1}\log(2(j+1)^2+2)\\
 &\leq\kappa_E(k_s+1)^2(J-j_0)\log(2J^2+2)\\
 &=O(J\log(J+2)).
 \end{aligned}
\]
Summing \eqref{eq:blockestimate} over $j_0\leq j<J$ therefore gives
\[
 \begin{aligned}
 K(B_{2J^2})
 &\geq K(B_{2j_0^2})
   +\left(1-\frac1{k_s}\right)\sum_{j=j_0}^{J-1}(4j+2)
   -\sum_{j=j_0}^{J-1}E_s(2(j+1)^2)\\
 &\geq K(B_{2j_0^2})
   +2\left(1-\frac1{k_s}\right)(J^2-j_0^2)
   -O(J\log(J+2))\\
 &=2\left(1-\frac1{k_s}\right)J^2-O(J\log(J+2)).
 \end{aligned}
\]
For $2J^2\leq n<2(J+1)^2$, the endpoint $2J^2$ is recursively
determined by $n$, so truncation gives
$K(B_n)\geq K(B_{2J^2})-O(1)$.
Since $n-2J^2=O(\sqrt n)$, it follows that
$\dimH(B)\geq1-\frac1{k_s}
 \geq1-\frac1{k_*}>\dimH(A)$,
a contradiction.
\end{proof}

\begin{proof}[Proof of Theorem~\ref{thm:main}]
Lemma~\ref{lem:domination} gives $A\leq_K B$.
For $N=2J^2$, summing \eqref{eq:blockestimate} over $0\leq j<J$, we obtain
\[
 \begin{aligned}
 K(B_N)
 &=K(B_0)+\sum_{j<J}
   \bigl(K(B_{2(j+1)^2})-K(B_{2j^2})\bigr)\\
 &\geq K(B_0)+\sum_{j<J}
   \left(1-\frac1{k_{s_j}}\right)(4j+2)
   -\sum_{j<J}E_{s_j}(2(j+1)^2)\\
 &\geq N-\sum_{j<J}\frac{4j+2}{k_{s_j}}
   -\sum_{j<J}E_{s_j}(2(j+1)^2),
 \end{aligned}
\]
where the last step uses $K(B_0)\geq0$.
By Lemma~\ref{lem:infinite}, $k_{s_j}\to\infty$, so the first
error sum is $o(J^2)$: for every $\varepsilon>0$, its tail is at
most $\varepsilon\sum_{j<J}(4j+2)=2\varepsilon J^2$.
Since each stage makes at most one repair, $s_j\leq j$, and hence
$k_{s_j}+1=O(\log(j+2))$. Thus the definition of $E_s(M)$ gives
\[
 E_{s_j}(2(j+1)^2)
 =\kappa_E(k_{s_j}+1)^2\log(2(j+1)^2+2)
 =O(\log^3(j+2)),
\]
with an implied constant independent of $j$. Summing over $j<J$, we obtain
\[
 \sum_{j<J}E_{s_j}(2(j+1)^2)
 =O(J\log^3(J+2))=o(J^2).
\]
Consequently, $K(B_N)\geq N-o(N)$.

For arbitrary large $n$, put $N=2\lfloor\sqrt{n/2}\rfloor^2$.
This endpoint is recursively determined by the input length $n$,
so $K(B_n)\geq K(B_N)-O(1)$ and $n-N=O(\sqrt n)$.
Thus $\dimH(B)=1$.

It remains to verify $A\leq_T B$. Fix the finite construction
parameters in the oracle program. Starting in state $(0,0)$, use
$B$ to scan the control positions for the first marker completion.
Read the four recorded integers $a$, $K(B_a)$, $L$, and $|p|$,
and then read $p=B[a+L:a+L+|p|)$. Running $p$ with condition
$(B_a,K(B_a))$ returns the exact word $A_{M+H_{s+1}}$; the needed
complexity value is supplied by the record. Compute
$M=2(\lfloor\sqrt{a/2}\rfloor+1)^2$, jump to $M$, and repeat in
state $(M,s+1)$. This procedure does not require the initial stored
word, since every recorded program returns an entire source prefix.
By Lemma~\ref{lem:infinite}, it recovers prefixes of unbounded
length, and hence computes every bit of $A$. The case $\dimH(A)=1$
was already handled by taking $B=A$.
\end{proof}

\section{Consequences and further questions}\label{sec:question}

We derive statistical and computational consequences of
Theorem~\ref{thm:main}, and then examine the question of $1$-random
$K$-upper bounds. We show that this question can be reduced to
generic sequences and that the randomness requirement cannot be
strengthened to weak Demuth randomness. Finally, we discuss analogous
questions for plain Kolmogorov complexity.

\subsection{Consequences of dimension one}

We first observe that Theorem~\ref{thm:main} yields $K$-upper bounds
satisfying the strong law of large numbers.
We say that $X\in\Cantor$ \emph{satisfies the strong law of large numbers} if
$\lim_{n\to\infty}\frac1n\sum_{i<n}X(i)=\frac12$.

\begin{corollary}\label{cor:strong-law}
For every $A\in\Cantor$, there is a $B\in\Cantor$ with $A\leq_K B$
that satisfies the strong law of large numbers.
\end{corollary}
\begin{proof}
Choose $B$ with $A\leq_K B$ and $\dimH(B)=1$ by
Theorem~\ref{thm:main}. Then there is a $1$-random
$R$ such that $\lim_{n\to\infty}\frac{|\{i<n:B(i)\neq R(i)\}|}{n}=0$
by \cite[Theorem~1.7]{GMSW}.
Every $1$-random sequence satisfies the strong law of large
numbers; see~\cite{DH,Nies}. Hence
\[
 \left|\frac1n\sum_{i<n}B(i)-\frac12\right|
 \leq\frac{|\{i<n:B(i)\neq R(i)\}|}{n}
 +\left|\frac1n\sum_{i<n}R(i)-\frac12\right|
 \to 0.
\]
\end{proof}

Combining the main theorem with Kurtz cofinality and Liu's theorem
gives a strict upper bound with both properties.

\begin{corollary}\label{cor:joint}
For every $A\in\Cantor$ there is a Kurtz random $R\in\Cantor$ such that $A<_K R$ and $\dimH(R)=1$.
\end{corollary}
\begin{proof}
By Liu's no-maximal-degree theorem~\cite[Theorem~1.1]{Liu}, choose $Y$
with $A<_K Y$. Apply Theorem~\ref{thm:main} to obtain $B$ with
$Y\leq_K B$ and $\dimH(B)=1$, and then apply Theorem~\ref{kz:cofinal}
to obtain a Kurtz random $R$ with $B\leq_K R$.
\end{proof}

A total function $f\colon\N\to\N$ is \emph{diagonally nonrecursive}
(DNR) if $f(e)\neq\varphi_e(e)$ whenever $\varphi_e(e)\downarrow$,
where $(\varphi_e)_{e\in\N}$ is a fixed standard enumeration of the
partial recursive functions. A set is of \emph{DNR degree} if it
computes a DNR function.

\begin{corollary}\label{cor:dnc}
For every
$A$, there is $B$ of DNR degree such that
$A\leq_K B$.
\end{corollary}
\begin{proof}
Choose $B$ with $A\leq_K B$ and $\dimH(B)=1$ by
Theorem~\ref{thm:main}.
Since $K(\sigma)\leq C(\sigma)+O(\log(|\sigma|+2))$, there is a
constant $d$ such that
$C(B\restr n)\geq\lfloor n/2\rfloor-d$ for every $n$.
Thus $B$ is \emph{complex}: the function $n\mapsto C(B\restr n)$
has a recursive nondecreasing unbounded lower bound.
By Kjos-Hanssen, Merkle and Stephan~\cite[Theorem~2.3]{KMS},
$B$ computes a DNR function.
\end{proof}

A sequence $(y_n)$ of real numbers is
\emph{uniformly distributed modulo one} if $$\lim_{N\to\infty}
 \frac{|\{n<N:\{y_n\}\in I\}|}{N}=\lambda(I)$$ for every interval
$I\subseteq[0,1]$,
where $\{y\}=y-\lfloor y\rfloor$ and $\lambda$ here denotes Lebesgue
measure. Following Avigad~\cite{Avigad}, a real $x$ is \emph{UD random}
if $(a_nx)_n$ is uniformly distributed modulo one for every recursive
sequence $(a_n)$ of pairwise distinct integers.

Recall that a \emph{Schnorr test} is a uniformly effectively open
sequence $(G_k)_{k\geq1}$ with $\lambda(G_k)\leq2^{-k}$ and with
$\lambda(G_k)$ uniformly recursive in $k$. A real passes the test if
it lies outside $\bigcap_kG_k$, and is \emph{Schnorr random} if it
passes every Schnorr test. Avigad~\cite[Theorem~2.1]{Avigad} proved
that every Schnorr random real is UD random. In particular, every
$1$-random real is UD random.

Effective Hausdorff dimension one and UD randomness are incomparable.
To see that dimension one does not imply UD randomness, set
$X(2^j)=0$ for every $j\geq0$ and fill the remaining positions, in
order, with the bits of a $1$-random sequence $R$.
Deleting the prescribed zero positions from $X_n$ recovers a prefix
of $R$ of length $n-O(\log(n+2))$. Hence $\dimH(X)=1$.
However, $\{2^{2^j}x_X\}<1/2$ for every $j$, so the recursive
sequence $a_j=2^{2^j}$ witnesses that $x_X$ is not UD random.
Conversely, by \cite[Theorem~7.4.8]{Nies}, there is a computably random sequence $Y$ such that
$K(Y_n)=O(\log(n+2))$.
Thus $\dimH(Y)=0$, whereas $Y$ is Schnorr
random and $x_Y$ is therefore UD random.

Nevertheless, effective dimension one implies uniform distribution
along every recursive sequence of pairwise distinct integers with
polynomial growth. The relevant classical result concerns the
Hausdorff dimension of the exceptional set.
Erd\H{o}s and Taylor~\cite[Theorem~13]{ET}
proved that, for an increasing integer sequence $a_n=O(n^d)$ with
$d>1$, the set of reals $x$ for which $(a_nx)_n$ is not uniformly
distributed modulo one has Hausdorff dimension at most $1-1/d$.
We make the covering argument effective: if $(a_n)$ is a recursive
sequence of pairwise distinct integers satisfying $|a_n|=O(n^d)$
for an integer $d\geq1$, then every exceptional
$x_X$ satisfies $\dimH(X)\leq1-1/d$.

\begin{proposition}\label{prop:uniform-distribution}
Let $X\in\Cantor$ satisfy $\dimH(X)=1$. Then $(a_nx_X)_n$ is
uniformly distributed modulo one for every recursive sequence
$(a_n)$ of pairwise distinct integers satisfying $|a_n|=O(n^d)$
for some integer $d\geq1$.
\end{proposition}
\begin{proof}
Fix a recursive sequence $(a_n)$ of pairwise distinct integers satisfying
$|a_n|=O(n^d)$ for an integer $d\geq1$, and choose an integer $L\geq1$
such that $\max_{n<N}|a_n|\leq LN^d$ for all $N\geq1$.
Suppose that $(a_nx_X)_n$ is not uniformly distributed
modulo one. By Weyl's criterion~\cite{Weyl},  there are an integer $h\ne0$ and
a rational $\delta>0$ such that
$$S_N(t)=\frac1N\sum_{n<N}e^{2\pi iha_nt}$$ satisfies $|S_N(x_X)|\geq4\delta$
for infinitely many $N$. Orthogonality of the distinct integer
frequencies $ha_n$ gives $\int_0^1|S_N(t)|^2\,dt=\frac1N$, and hence 
\[
 \lambda\{t\in[0,1]:|S_N(t)|>\delta\}
 \leq\frac1{N\delta^2}.
\]
We turn this estimate into an effective cover by short dyadic intervals.

Put $q=1+\delta/4$ and $N_j=\lceil q^j\rceil$ for $j\geq1$.
For $M\geq N\geq1$, splitting the sum at $N$ gives
\[
 S_M(t)-S_N(t)
 =\left(\frac1M-\frac1N\right)\sum_{n<N}e^{2\pi iha_nt}
   +\frac1M\sum_{N\leq n<M}e^{2\pi iha_nt}.
\]
Since each exponential has modulus one, the triangle inequality yields
\[
 |S_M(t)-S_N(t)|
 \leq\left(\frac1N-\frac1M\right)N+\frac{M-N}{M}
 =\frac{2(M-N)}M.
\]
As $N_j/N_{j-1}\to q$, whenever $N_{j-1}<N\leq N_j$ and $N$
is sufficiently large, we have
\[
 |S_{N_j}(t)-S_N(t)|
 \leq 2\left(1-\frac{N_{j-1}}{N_j}\right)<\delta.
\]
Hence
$|S_{N_j}(x_X)|\geq3\delta$ for infinitely many $j$.
The polynomial growth assumption gives
\[
 \|S_N'\|_\infty
 \leq\frac{2\pi|h|}{N}\sum_{n<N}|a_n|
 \leq 2\pi|h|LN^d.
\]
Fix a sufficiently large nonnegative integer $c$ and put
$m_j=d\log N_j+c$, so that
$2\pi|h|LN_j^d2^{-m_j}\leq\delta/4$ for every $j$.
For each $\sigma\in2^{m_j}$, the dyadic interval and its midpoint are $I_\sigma=[0.\sigma,0.\sigma+2^{-m_j}]$ and $t_\sigma=0.\sigma+2^{-m_j-1}$.
Compute a rational $r_{j,\sigma}$ satisfying $\bigl|r_{j,\sigma}-|S_{N_j}(t_\sigma)|\bigr|\leq\delta/4$,
and set
\[
 F_j=\{\sigma\in2^{m_j}:r_{j,\sigma}>2\delta\}.
\]
These finite sets are uniformly recursive. Since
\[
 |S_{N_j}(t)-S_{N_j}(t_\sigma)|\leq\delta/4
 \qquad(t\in I_\sigma),
\]
we have
\[
 \begin{gathered}
 |S_{N_j}(x_X)|\geq3\delta\ \Longrightarrow\ X_{m_j}\in F_j,\\
 \bigcup_{\sigma\in F_j}I_\sigma
 \subseteq\{t\in[0,1]:|S_{N_j}(t)|>\delta\}.
 \end{gathered}
\] The measure bound therefore yields $|F_j|2^{-m_j}\leq\frac1{\delta^2N_j}$.
For each of the infinitely many $j$ with $X_{m_j}\in F_j$,
describing $j$ and the index of $X_{m_j}$ in $F_j$ gives
\[
 \begin{aligned}
 K(X_{m_j})
 \leq K(j)+\log|F_j|+O(1)
 \leq m_j-\log N_j+O(\log(j+2)).
 \end{aligned}
\]
Since $\log N_j=j\log_2q+O(1)$ and
$m_j=d\log N_j+O(1)$, dividing the preceding inequality by $m_j$ gives
$\dimH(X)\leq1-1/d<1$, a contradiction.
\end{proof}

Combining Proposition~\ref{prop:uniform-distribution} with
Theorem~\ref{thm:main} yields the corresponding $K$-upper bounds.

\begin{corollary}\label{cor:uniform-distribution}
For every $A\in\Cantor$, there is a $B\in\Cantor$ with $A\leq_K B$
such that $(a_nx_B)_n$ is uniformly distributed modulo one for every
recursive sequence $(a_n)$ of pairwise distinct integers satisfying
$|a_n|=O(n^d)$ for some integer $d\geq1$.
\end{corollary}

\subsection{The question of a 1-random upper bound}
We return to the problem in~\eqref{eq:random-upper-question}.
The existence of random upper bounds can be tested on generic
sequences. For $S\subseteq\Strings$ and $G\in\Cantor$, $G$ meets $S$
if $(\exists n)\,G_n\in S$, and $G$ avoids $S$ if
$(\exists n)(\forall\tau\succeq G_n)\,\tau\notin S$.
The set $S$ is dense if every finite string has an extension in $S$.
A sequence is \emph{$1$-generic} if it meets or avoids every r.e.\ set
of finite strings, and \emph{weakly $2$-generic} if it meets every
dense set of finite strings that is r.e.\ relative to $\emptyset'$.
A sequence $G$ is \emph{low} if $G'\equiv_T\emptyset'$.

There is a $1$-generic sequence with a $1$-random $K$-upper bound.
Indeed, choose a nonrecursive r.e.\ $K$-trivial set $C$.
Every nonrecursive r.e.\ set computes a $1$-generic sequence, so
there is a $1$-generic $G\leq_T C$. Then $G$ is $K$-trivial.
The following proposition shows that the analogous existence statement
for weakly $2$-generic sequences is equivalent to
\eqref{eq:random-upper-question}.

\begin{proposition}\label{prop:generic-reduction}\begin{enumerate}
    \item If there is a weakly $2$-generic sequence with a $1$-random
$K$-upper bound, then every $K$-degree has a $1$-random upper bound.

\item If some $K$-degree has no $1$-random upper bound,
then there is a low $1$-generic sequence with no $1$-random
$K$-upper bound.
\end{enumerate}
\end{proposition}
\begin{proof}
Let $\mathcal U$ be the class of sequences having a $1$-random
$K$-upper bound. For $c,d\in\N$, let $\mathcal P_{c,d}$ consist of the
sequences $X$ for which some $R\in\Cantor$ satisfies $K(R_n)\geq n-d$ and $K(X_n)\leq K(R_n)+c$ for every $n$.
By definition, $\mathcal U=\bigcup_{c,d\in\N}\mathcal P_{c,d}$ and the classes $\mathcal P_{c,d}$ are uniformly
$\Pi^0_1(\emptyset')$. 
The class $\mathcal U$ is invariant under finite changes. If
$\mathcal U\neq\Cantor$, its complement is therefore dense. Let $T_{c,d}$ be the $\emptyset'$-recursive tree such that $[T_{c,d}]=\mathcal P_{c,d}$. 
Then the uniformly $\emptyset'$-recursive sets of strings $D_{c,d}=\Strings\setminus T_{c,d}$
are dense. Every weakly $2$-generic meets every $D_{c,d}$ and hence belongs
to none of the classes $\mathcal P_{c,d}$. This proves (1).

Under the same assumption $\mathcal U\neq\Cantor$, fix an effective
enumeration $(W_e)_e$ of the r.e.\ sets of finite strings and a recursive
enumeration $((c_e,d_e))_e$ of $\N^2$. Using $\emptyset'$, construct
an increasing sequence $(\sigma_e)$, starting with the empty word.
At step $e$, extend $\sigma_e$ into $W_e$ if such an extension exists,
and otherwise keep $\sigma_e$. Then extend into $D_{c_e,d_e}$,
increasing the length, to obtain $\sigma_{e+1}$. The union
$G=\bigcup_e\sigma_e$ is $1$-generic, recursive in $\emptyset'$,
and belongs to none of the classes $\mathcal P_{c,d}$.
Since every $\Delta^0_2$ $1$-generic sequence is low, this proves (2).
\end{proof}

The randomness requirement in \eqref{eq:random-upper-question} cannot, however,
be strengthened to weak Demuth randomness. This obstruction already
occurs for a non-1-random $K$-degree.

Fix an effective enumeration $(\mathcal W_e)_{e\in\N}$ of the r.e.\ open
subsets of $\Cantor$. A function $g:\N\to\N$ is
\emph{$\omega$-recursively approximable} if it has a recursive
approximation $(g_s)_{s\in\N}$ and a total recursive function $b$ such that $g(n)=\lim_sg_s(n)$ and $|\{s:g_{s+1}(n)\neq g_s(n)\}|\leq b(n)$ for all $n$.
A \emph{Demuth test} is a sequence
$(\mathcal V_n)_{n\in\N}$ of the form
$\mathcal V_n=\mathcal W_{g(n)}$, where $g$ is
$\omega$-recursively approximable and $\lambda(\mathcal V_n)\leq2^{-n}$.
A sequence $X$ is \emph{weakly Demuth random} if $X\notin\bigcap_n\mathcal V_n$ for every Demuth test $(\mathcal V_n)_{n\in\N}$.
This definition is equivalent to restricting to tests with nested
final components: replace $\mathcal V_n$ by
$\bigcap_{i\leq n}\mathcal V_i$.
See Ku\v{c}era and Nies~\cite[Section~2.2]{KN11}.

Fix a universal prefix-free machine $U$, Chaitin's halting
probability is
$\Omega=\sum_{p\in\dom U}2^{-|p|}$,
identified with its binary expansion. It is $1$-random and
\emph{left-r.e.}, meaning that it is the limit of a recursive
nondecreasing sequence of rational numbers.
We continue to write $\Omega_N=\Omega\restr N$.

\begin{proposition}\label{prop:weak-demuth-obstruction}
There is a non-$1$-random sequence $A\in\Cantor$ with no weakly
Demuth random $K$-upper bound. Moreover, $A$ can be chosen so that
$A<_K\Omega$ and $\dimH(A)=1$.
\end{proposition}
\begin{proof}
Let $D=\{4^j:j\geq1\}$ and $d(N)=|D\cap[0,N)|$, and define
\[
 A(i)=\begin{cases}
 0,&i\in D,\\
 \Omega(i),&i\notin D.
 \end{cases}
\]
For some constant $b$,
\[
 K(\Omega_N)-d(N)-b\leq K(A_N)\leq K(\Omega_N)+b
 \qquad(N\in\N).
\]
Since $A$ is zero on the infinite recursive set $D$, it is not
$1$-random, and hence $A<_K\Omega$. Also,
$d(N)=O(\log(N+2))$, so $\dimH(A)=1$.

By the ample excess lemma~\cite{MY08},
\[
 L:=\sum_{N\in\N}2^{N-K(\Omega_N)}<\infty.
\]
For $n\geq1$ and $\tau\in2^{4n}$, let $m(\tau)$ be the integer
with binary representation $1\tau$, and put
\[
 \mathcal U_n=\bigcup\bigl\{[\tau]:\tau\in2^{4n},
 K(A_{m(\tau)})\leq m(\tau)+n\bigr\}.
\]
The map $m$ sends $2^{4n}$ bijectively onto
$[2^{4n},2^{4n+1})\cap\N$. The cylinders defining $\mathcal U_n$
are pairwise disjoint and each has measure $2^{-4n}$.
For every included cylinder $[\tau]$, putting $N=m(\tau)$ gives
$K(A_N)\leq N+n$, and hence
\[
 2^{N-K(A_N)}\geq2^{-n},\qquad
 2^{-4n}\leq2^{-3n}2^{N-K(A_N)}.
\]
Summing over these cylinders and enlarging the sum to the whole
interval gives the first inequality below.
On this interval, $d(N)\leq2n$ and
hence $2^{N-K(A_N)}\leq2^{b+2n}2^{N-K(\Omega_N)}$.
Therefore
\[
 \begin{aligned}
 \lambda(\mathcal U_n)
 \leq2^{-3n}\sum_{N=2^{4n}}^{2^{4n+1}-1}2^{N-K(A_N)}
 \leq2^{b-n}\sum_{N=2^{4n}}^{2^{4n+1}-1}2^{N-K(\Omega_N)}
 \leq2^bL\,2^{-n}.
 \end{aligned}
\]

Suppose $A\leq_KX$, and let $m=m(X_{4n})$. A prefix-free machine
first reads a shortest description of $n$, then the $4n$ bits of
$X_{4n}$, which determine $m$, and finally the remaining $m-4n$ bits.
Thus
\[
 K(A_m)\leq K(X_m)+O(1)
 \leq m+K(n)+O(1)\leq m+n
\]
for all sufficiently large $n$. Hence $X\in\mathcal U_n$ for all
sufficiently large $n$.

Since $\Omega$ is left-r.e., its length-$N$ prefixes have recursive
approximations with at most $2^N-1$ changes. Overwriting the positions
in $D$ gives approximations $A_{N,s}$ with the same bound. Use these
and the usual decreasing approximation $K_s$ to obtain finite clopen
approximations $\mathcal U_{n,s}$ to $\mathcal U_n$. While $A_{N,s}$
is unchanged, the inequality $K_s(A_{N,s})\leq N+n$ changes its truth
value at most once. Thus the finite sets of cylinders defining
$\mathcal U_{n,s}$ have a recursive bound on their number of changes,
uniformly in $n$, and so do their canonical r.e.\ open indices.
Choose a fixed $q\geq1$ with $2^bL\leq2^q$. Then
$(\mathcal U_{n+q})_{n\in\N}$ is a Demuth test. Every $K$-upper bound
of $A$ belongs to all components of some tail of this test. Since
every tail is again a Demuth test, no such upper bound is weakly
Demuth random.
\end{proof}

\subsection{Plain complexity}\label{sec:plain}
Replacing prefix-free complexity by plain complexity gives the
reducibility
\[
 A\leq_C B\quad\Longleftrightarrow\quad
 (\exists c)(\forall n)\ C(A_n)\leq C(B_n)+c.
\]
The equivalence classes under mutual reducibility are the
\emph{$C$-degrees}. We write $A<_C B$ when $A\leq_C B$ but
$B\nleq_C A$. The results above suggest several analogues for this
degree structure.

\begin{question}\label{qu:plain-kurtz}
Does every $C$-degree have a Kurtz random upper bound?
\end{question}

The difficulty in adapting Theorem~\ref{kz:cofinal}
is removing the auxiliary nonrecursive oracle. That proof uses $K$-triviality, whereas a set $W$ satisfying
$C(W_n)\leq C(n)+O(1)$ must be recursive by Chaitin \cite{Chaitin76}.

We also conjecture that every $C$-degree has an upper bound of effective
Hausdorff dimension one.

\section*{acknowledgments}
The author was supported by Focused Research Grants of Natural Science Foundation of Jiangsu Province No. BK20243060.  The author thanks Professor Lu Liu and Liang Yu for helpful discussions.
The author used OpenAI’s ChatGPT 6 to assist with developing the proofs. All proofs in this paper have been independently verified by the author without the use of AI tools. The author takes full responsibility for the final content of the paper.

\end{document}